\documentclass[11pt,a4paper]{article}
\usepackage{amsfonts,amsmath,amstext,amssymb}
\usepackage{dsfont}
\usepackage{theorem}    
\usepackage{a4wide}
\usepackage{verbatim}

\newtheorem{lema}{Lemma}
\newtheorem{prop}[lema]{\bf Proposition}
\newtheorem{teor}[lema]{\bf Theorem}
\newtheorem{coro}[lema]{\bf Corollary}
\newtheorem{ejem}[lema]{\bf Example}
\newtheorem{rema}[lema]{\bf Remark}

\theorembodyfont{\rmfamily}

\renewcommand{\sup}{\mathrm{Sup}}

\title{Uniqueness and non-existence of minimal submanifolds}

\author{
Rafael M. Rubio and  Juan J. Salamanca\\[6mm]
Departamento de Matem\'aticas, Campus de Rabanales, \\[0.5mm]
Universidad de C\'ordoba, 14071 C\'ordoba, Spain,\\[0.5mm]
E-mails\textup{:\texttt{\;rmrubio@uco.es},\,\,\texttt{jjsalamanca@uco.es}}}

\date{}

\begin{document}

\maketitle

\thispagestyle{empty}

\begin{abstract}
We provide uniqueness results for compact minimal submanifolds in a large
class of Riemannian manifolds of arbitrary dimension. In the case compact and Cartan-Hadamard
manifolds we obtain general results for these submanifolds. Several applications to Geometric
Analysis are also showed.

\vspace*{2mm}

\end{abstract}



\vspace*{-2mm}

\section{Introduction}

The importance of minimal submanifolds (and, in particular, minimal surfaces) is very well-known. 
Since historical reasons, the problem of minimal hypersurfaces was firstly studied as graphs in $\mathbb{R}^{n+1}$. 
That is, given a function $u \in C^\infty (\Omega)$,
$\Omega$ an open domain in $\mathbb{R}^n$, the graph $\Sigma_{u}= \left\{ (u(p),p): p\in \Omega \right\}$ in the Euclidean
space $\mathbb{R}^{n+1}$
defines a minimal hypersurface. It can be shown that $u$ defines such a minimal graph if and only if $u$ satisfies the following quasi-linear elliptic PDE,
\begin{equation}\label{euclidean}
\mathrm{div}\left( \frac{Du}{\sqrt{1+|Du|^2}} \right) = 0 \, .
\end{equation} This equation has been widely studied and even nowadays much researchers pay attention to it. From an analytical
point of view, it is the Euler-Lagrange equation of a classical variational problem. For each $u \in C^\infty (\Omega)$,
$\Omega$ an open domain in $\mathbb{R}^n$, the volume element of the induced metric from $\mathbb{R}^{n+1}$ is represented
by the $n$-form $\sqrt{1+|Du|^2} dV$ on the graph $\Sigma_{u}$, where $dV$ is the
volume form of $\Omega$. The critical points
of the $n$-volume functional $u \mapsto \int \sqrt{1+|Du|^2} dV$ are given by the equation (\ref{euclidean}).
In 1914, S. Bernstein \cite{Bernstein}, amended latter by E. Hopf in 1950 \cite{Hopf}, proved his well-known
uniqueness theorem for $n=2$,
\begin{quote}{\it
The only entire solutions to the minimal surface equation in $\mathbb{R}^3$ are the affine functions
$$
u(x,y)=ax+by+c\, ,
$$ where $a,b,c\in \mathbb{R}$.
} 
\end{quote} 
In terms of PDEs, Bernstein proved a general Liouville type result,
\begin{quote}{\it
Any bounded solution $u \in C^\infty (\mathbb{R}^2)$ of the PDE
$$
A \, u_{_{xx}} + 2 B \, u_{_{xy}} + C \, u_{_{yy}} =0 \, ,
$$ where $A,B,C \in C^\infty (\mathbb{R}^2)$ such that $AC-B^2 >0$, must be constant.
}
\end{quote} Then, a lot of work has been made in order to extend the classical Bernstein result to higher dimensions (see \cite{Osserman}
for a survey until 1984). A notable progress was made by J. Moser \cite{Moser} in 1961, who obtained the so-called Moser-Bernstein theorem,
\begin{quote}{\it
The only entire solutions $u$ to the minimal surface equation in $\mathbb{R}^{n+1}$ such that $|Du|\leq C$, for some $C \in \mathbb{R}^+$, are the
affine functions
$$
u(x_1, \ldots, x_n) = a_1 x_1 + \ldots + a_n x_n + c \, ,
$$ where $a_i, c \in \mathbb{R}$, $1\leq i \leq n$ and $\sum_{i=1}^n a_i^2 \leq C^2$.
}
\end{quote} In 1951 L. Bers, \cite{Bers}, proved that a solution $u$ of the minimal surface equation in $\mathbb{R}^3$ defined
on the exterior of a closed disc in $\mathbb{R}^2$ has bounded $|Du|$. Hence, the Moser-Bernstein theorem for $n=2$ and the Bers
result provide another proof for the classical Bernstein theorem.
In 1968, J. Simons \cite{Simons}, together with other results of E. De Giorgi \cite{Giorgi} and W.H. Fleming \cite{Fleming} yield a proof
of the Bernstein theorem for $n \leq 7$. Furthermore, it was found a counterexample $u \in C^\infty (\mathbb{R}^n)$ for each $n \geq 8$.

Then, much research has been made in order to characterize minimal submanifolds. For example, H. Rosenberg study
minimal surfaces in the product of $\mathbb{R}$ and a Riemannian surface in \cite{Rosenberg}. Minimal surfaces are
studied in warped product manifolds in several papers (see, for instance, ).

It is known that $\mathbb{R}^n$ does not admit any compact minimal submanifold. However, this fact does not occur in $\mathbb{S}^n$. Hence, it 
is natural to consider the problem of obtaining characterization results for minimal submanifolds in a large class of Riemannian manifolds. We 
require that the structure of the Riemannian manifold splits into an open interval of the real line and a Riemannian manifold. In the end, we will see
that this topological assumption can be avoided in some relevant cases (see Corollary \ref{corCH} for instance).

Consider
a smooth $1$-parametrized family of Riemannian metrics $(F^{n+m},g_{_t})$, $t\in I \subseteq \mathbb{R}$, on a differential manifold $F$, and a
positive function $\beta \in C^\infty (I)$. The product 
manifold $I\times F$ can be endowed with the following metric
\begin{equation}
 \overline{g}=\beta \, \pi_{_I}^*(dt^2) +\pi_{_F}^*(g_{_t}) \quad \textrm{(} \, \beta \,  dt^2 +g_{_t} \quad  \textrm{in short)} \, ,
\end{equation} where $\pi_{_I}$ and $\pi_{_F}$ denote the canonical projections onto $I$ and $F$ respectively. 

Observe that a suitable open normal neighbourhood of an arbitrary Riemannian manifold lies in this family (in this case, consider $\beta=1$ and $t$ as
the geodesic distance to a fixed point of the neighbourhood). In particular, removing a point of a
simply-connected complete Cartan-Hadamard manifold, we have that the resulting manifold possesses this structure.

Moreover, this kind of Riemannian manifolds generalizes properly to the important class of warped product Riemannian manifolds (see, for instance \cite{ON}).

We focus our attention on the case in which the Riemannian manifold has an isotropic behaviour associated to the $t$ coordinate. Given any 
compact subset, we desire that, by the flux along $\partial_t$, its volume does not increase or decrease. To make clear this
idea, we may introduce the following notion: a Riemannian manifold of the form $(I\times F,\beta dt^2+g_{_t})$
is \emph{non-shrinking} (resp. \emph{expanding}) \emph{throughtout} $\partial_t$ if
$$
\partial_t \beta \geq 0 \, \, \, \quad {\rm and} \quad (\mathcal{L}_t g_{_t})(X,X) \geq 0 \, \, \, ({\rm resp.} > 0) \, ,
$$  for any $X\in \mathfrak{X}(F)$. For the non-shrinking case, this definition is equivalent to $\mathcal{L}_{\partial_t}\overline{g}$ to be a definite non-negative
tensor field.  From now on, $\partial_t g_{_t}$ will mean $\mathcal{L}_{\partial_t} g_{_t}$. Taking $p,q \in F$ close enough, the distance in $F$ between $p_t:=(t,p)$ and $q_t:=(t,q)$ measured at some value of $t$
is given by $d_{_t}(p_t,q_t)$, where $d$ is the induced Riemannian distance. Then, the Riemannian manifold is non-shrinking (resp. expanding) 
throughtout $\partial_t$ if
$d_{_t}(p_t,q_t)$ is a non-decreasing (resp. increasing) function of $t$. Dually, we will say that a Riemannian manifold is \emph{non-expanding}
(resp. \emph{contracting}) throughout $\partial_t$ if it is non-shrinking (resp. expanding) throughout $-\partial_t$. The geometrical
interpretations also hold with their respectives changes. We may unify these two notions with the following one. We 
will say that the manifold $(I\times F, \beta \, dt^2+g_{_t})$ is \emph{monotone} (resp. \emph{strictly monotone}) if it is non-shrinking 
or non-expanding (resp. expanding or contracting) throughtout $\partial_t$.

We will prove (see Proposition \ref{local}) that, locally, any Riemannian manifold is locally expanding throughtout certain vector field. Moreover,
it can be be shown that 
any simply-connected complete Cartan-Hadamard manifold removing a point is also expanding throughout a vector field (Proposition \ref{CH}).

This paper is organized as follows. In Section 3 we show several characterization results for the class of Riemannian manifolds introduced.
The first of them is Theorem \ref{T1},

\begin{quote}
{\it
In a monotone Riemannian manifold $(I\times F,\beta \, dt^2 + g_{_t})$, every compact minimal submanifold
must be contained in a level hypersurface $t= const$.
Moreover, in the case of codimension higher than one, $S$ is a minimal submanifold of $(F,g_{_{t_0}})$, for some
$t_{_0}\in I$.
}
\end{quote}

We analyze several consequences. Then, in Section 4 we study the behavior of a distinguished function, what lead us to 
some non-existence results. In Section 5 we show how our techniques lead to relate symmetries of the Riemannian manifold 
with minimal submanifolds, Theorem \ref{TSimetrias},
\begin{quote} 
{\it 
Let $M$ be a complete simply-connected Riemannian manifold which admits an irrotational nowhere zero Killing vector field $K$. Every minimal compact submanifold must be contained in a leaf of the foliation $K^\bot$.
}
\end{quote} Then, Section 6 is devoted to present some uniqueness results for a wide family of PDEs, leading to solve new Bernstein type problems.
 In Section 7 we provide some results for certain Dirichlet problems. After that, in the last section we focus our attention on two families of Riemannian manifolds: the compact and the Cartan-Hadamard. For the compact case, we obtain that any
compact Riemannian manifold admits a positive number $\delta$ such that any compact minimal submanifold cannot be contained in a open geodesic ball of radius
$\delta$. We denote $\mathrm{diam}(M,g):=
\max \left\{ d(p,q) , p,q\in M \right\}$.
A normalization shows that the quantity $\overline{\delta}(M,g):= \sup \delta :$ there exits no compact minimal submanifold contained
in a geodesic ball of radius $\delta /\mathrm{diam}(M,g) $ is upper bounded by $1$. In fact, we prove that $\overline{\delta}(M,g)\in (0,1)$. We get that $\overline{\delta}(\mathbb{S}^n,g_{\mathbb{S}^n})=1/2$ (see Corollary \ref{curvaturaconstante}). We inquire if this is a characterization result, that is, is the round sphere the only compact Riemannian manifold $(M,g)$ such that $\overline{\delta}(M,g)=1/2 $?

Finally, we consider the class of Cartan-Hadamard manifolds. We obtain that any Cartan-Hadamard manifold removing a point $p$ is expanding
throuthout the polar normal coordinate vector field centered at $p$. We get Theorem \ref{TeoCH} (compare with \cite[Corollary 2]{atsuji}),
\begin{quote} {\it
Let $(M,g)$ be a simply-connected complete Cartan-Hadamard Riemannian manifold. It admits no compact minimal submanifold.}
\end{quote} As a consequence, the simply-connected assumption can be assumed on the compact minimal submanifold. As a particular case, any
topological $n$-sphere cannot be immersed minimally
in a Cartan-Hadamard manifold. We end this paper with an analysis on the shape that a minimal submanifold can have, Theorem \ref{noextremos},

\begin{quote}{\it
Let $x:S \rightarrow M^n$ be a minimal submanifold in a complete Cartan-Hadamard manifold. The lift of the minimal immersion in the
universal Riemannian covering, $\tilde{x}:\tilde{S} \rightarrow \mathbb{R}^n$ cannot have a strict extremum point.}
\end{quote}
We also show how to construct Riemannian manifolds where this fact does not hold.

\section{Preliminares}

Consider the Riemannian manifold $(\overline{M}=I\times F,\overline{g}=\beta dt^2+g_{_t})$. Let $x:S^n \rightarrow I\times F$ be 
an $n$-dimensional submanifold. On $S$, take
the function $\tau:=\pi_{_I} \circ x$, where $\pi_{_I}$ is the projection onto $I$. It is not difficult to obtain that its gradient satisfies
\begin{equation}\label{gradiente}
\nabla \tau =\frac{1}{\beta} \partial_t^\top \, . 
\end{equation} On $S$, define the acute angle function $\theta$ ($\theta \in \left[0,\pi \right]$)  between $S$ and $\partial_t$, by
\begin{equation}
|\nabla \tau|^2 =: \frac{1}{\beta} \sin^2 \theta \, ,
\end{equation} where $\nabla$ denotes here the gradient on $S$. Equivalently, at each point of $S$, $\sin^2 \theta= |\frac{1}{\sqrt{\beta}} \partial_t^\top|^2$. Hence, this function is well-defined. Clearly, 
when $S$ is an hypersurface, $\theta$ is the angle between
the normal vector field and the unit vector field $\frac{1}{\sqrt{\beta}}\partial_t$, i.e., $\cos \theta = \overline{g}(N,\frac{1}{\sqrt{\beta}} \, \partial_t)$.

Note that $S$ is
contained in an hypersurface $t=const$ 
if and only if $\theta$ vanishes identically. Let us write
\begin{equation} \label{decompositiontt}
\partial_t^\top = \partial_t -\sum_{i=1}^m \overline{g}\left( N_i,\partial_t \right) N_i \, , 
\end{equation} where $\left\{ N_i \right\}_{i=1}^m$ is an orthonormal basis of $T_p^\bot S$, $p\in S$, and $m$ is the codimension of $S$. From (\ref{gradiente}),
\begin{equation} \label{coseno}
 \cos^2 \theta = \sum_{i=1}^m \frac{1}{\beta} \overline{g} \left( N_i,\partial_t  \right)^2 \, .
\end{equation}  

 In another setting, recall that the second fundamental form tensor, 
$\mathrm{II}:\mathfrak{X}(S) \times \mathfrak{X}(S) \rightarrow \mathfrak{X}^\bot (S)$ of $S$ is
$$
\mathrm{II}(X,Y)= \left( \overline{\nabla}_X Y \right)^\bot \, ,
$$ for any $X,Y\in \mathfrak{X}(S)$, where $\overline{\nabla}$ is the Levi-Civita 
connection of the ambient space. A contraction of this tensor produces the mean curvature vector field $H$. 
Namely, if $\left\{ E_i \right\}_{i=1}^n$ is an orthonormal basis of $T_p S$, $p\in S$,
$$
nH = \sum_{i=1}^n \mathrm{II}\left( E_i,E_i \right) \, .
$$ A submanifold is said to be \emph{minimal} provided that $H=0$.

\section{Main Results}

First, we obtain some useful formulae. We begin computing
the Laplacian of a distinguished function on minimal submanifolds.

Consider a Riemannian  submanifold $(S,g)$ of the Riemannian manifold $(\overline {M}=I\times F,\overline{g}=\beta dt^2+g_{_t})$. 
Let $\left\{ E_i \right\}_{i=1}^n$ and be a local frame of $S$ on an open set ${\cal U\subset} S$ and let  $\left\{ N_j \right\}_{j=1}^m$ be 
a local frame on ${\cal U\subset} S$ of the normal vector bundle of $S$ in $\overline{M}$. Standard computations, making use of (\ref{decompositiontt}), lead 
to the following expression of the Laplacian of $\tau$ in $(S,g)$,
\begin{equation}\label{L}
\Delta \tau = -\frac{1}{\beta^2} \partial_t^\top (\beta) + \frac{1}{\beta} \mathrm{\overline{div}}\left( \partial_t \right)
- \sum_{i=1}^m \frac{1}{\beta} \overline{g} \left( \overline{\nabla}_{N_i} \partial_t,N_i \right) 
- \sum_{i=1}^m \sum_{j=1}^n \overline{g}\left(N_i, \partial_t \right) \overline{g} \left( \overline{\nabla}_{E_j} N_i, E_j \right) \, ,
\end{equation}

\noindent where $ \mathrm{\overline{div}}$ and $\overline{\nabla}$ denote the divergence operator and the Levi-Civita connection of $(\overline{M},\overline{g})$ respectively.

It is not difficult to show
$$
\overline{g} \left( \overline{\nabla}_{E_j} N_i, E_j \right) = - \overline{g}\left( N_i, \mathrm{II}(E_j,E_j) \right) \, .
$$ Therefore,  the last addend of (\ref{L}) vanishes when   $S$ is minimal.

On the other hand, let us write
$N_i=\frac{1}{\beta} \overline{g}\left( N_i,\partial_t \right) \partial_t +N_i^F \, ,
$ where $\overline{g} \left( N_i^F ,\partial_t  \right)=0$. Then, each term of the form $\overline{g} \left( \overline{\nabla}_{N_i} \partial_t,N_i \right)$ can be
decomposed as follows. Let $p \in S$ be such that $N_i^F (p) \neq 0$. Take us new coordinates (reducing the size of ${\cal U}$, if  it is necessary)
$(\mathcal{U}, (t\equiv x_0,x_1,\ldots x_{n+m-1}) )$ around $p$ in $I\times F$ such that $N_i^F = \partial_{x_1}$ on $S \cap \mathcal{U}$.
From the definition of the Christoffel symbols for a coordinate system, we have
$$
\overline{g}\left( \overline{\nabla}_{\partial_t} \partial_t,N_i^F  \right) = -\frac{1}{2} \partial_{x_1} \overline{g}(\partial_t,\partial_t) 
$$ and $$
\overline{g}\left( \overline{\nabla}_{N_i^F} \partial_t,\partial_t  \right) =  \frac{1}{2} \partial_{x_1} \overline{g}(\partial_t,\partial_t)\, . 
$$ Moreover,
$$
\overline{g}\left( \overline{\nabla}_{N_i^F} \partial_t, N_i^F \right) =\frac{1}{2} \partial_t g_{_t} (N_i^F,N_i^F) \, .
$$ 

\noindent Note that the previous equations also hold at every point $p\in {\cal U}$ with $N_i^F(p)=0$. Hence, we arrive to
$$
\overline{g} \left( \overline{\nabla}_{N_i} \partial_t,N_i \right) = 
 \frac{\overline{g}(N_i,\partial_t)^2}{2 \beta} \partial_t \beta + \frac{1}{2} (\partial_t g_{_t} )(N_i^F,N_i^F) \, .
$$ 

Now, to analyze $\mathrm{\overline{div}}\left( \partial_t \right)$ a more suitable frame field is needed. At each $p=(t_{_0},q) \in I\times F$ take normal coordinates around $p$ such that  $\left\{ \frac{1}{\sqrt{\beta}} \partial_t (p), \partial_{x_{_i}}(p) \right\}_{i=1}^{n+m-1}$ is an orthonormal basis on the tangent 
space at $p$. The expression of the metric allows us to obtain
\begin{equation} \label{divt}
\mathrm{\overline{div}}(\partial_t) = \frac{1}{2 \beta} \partial_t \beta +\sum_i  \frac{1}{2} \partial_t g_{_{ii}} \, ,
\end{equation} where the $g_{_{jk}}$ are the components of the metric tensor $g_{_t}$ in this coordinate chart.
From previous equation, it is clear that if the Riemannian manifold is non-shrinking (resp. non-expanding)
throughout $\partial_t$, then $\mathrm{\overline{div}}(\partial_t)$ is a non-negative (resp. non-positive) function.

Therefore, taking into account all previous considerations, 
\begin{equation} \label{laplaciano}
 \Delta \tau = -\frac{1}{\beta^2} \partial_t^\top \beta +\frac{1}{2 \beta} \left[  
\left(1-\sum_j \overline{g}(N_j, \frac{1}{\sqrt{\beta}}\partial_t)^2 \right) \, \frac{\partial_t \beta}{\beta}  + \sum_i \partial_t g_{_{ii}} - (\partial_t g_{_t})(N_i^F, N_i^F)
 \right] \, .
\end{equation} Now, assuming that $S$ has dimension at least three, we can consider the 
following pointwise conformal metric on $S$, $\widetilde{g}= \beta^{(n-2)/2} g$. Endowed $S$ with this metric and making use of (\ref{coseno}), then
the $\widetilde{g}$-Laplacian of $\tau$ becomes
\begin{equation} \label{laplacianoconforme}
 \widetilde{\Delta} \tau = \frac{\beta^{-n/2}}{2} \left[  
\sin^2 \theta \, \frac{\partial_t \beta}{\beta}  + \sum_i \partial_t g_{_{ii}} - (\partial_t g_{_t})(N_i^F, N_i^F)
 \right] \, .
\end{equation}

\begin{rema} \label{dimensionmas}
{\rm Although this conformal change does not apply when the submanifold is $2$-dimensional, nevertheless, we 
can  build a $1$-dimensional extension in order to increase the dimension of the submanifold
and the Riemannian manifold. Indeed, if $x:S \rightarrow I \times F$ is a minimal $2$-dimensional 
isometric immersion, we consider the Riemannian manifold
$\left( (I\times F)\times \mathbb{S}^1,\overline{g}+ds^2\right)$, being $ds^2$ the standard metric of $\mathbb{S}^ 1$ and the  following
natural $3$-dimensional isometric immersion, $\hat{x}: S\times \mathbb{S}^1 \rightarrow \left( I\times F\right) \times \mathbb{S}^1$, with
$\hat{x}(p,s)=(x(p),s)$, for all $p\in S$ and $s\in \mathbb{S}^1$. 

Taking into account the natural identifications $T_{(p,\alpha)}(S\times \mathbb{S}^1)\equiv T_pS\oplus T_\alpha \mathbb{S}^1,$ $p\in S$, $\alpha\in \mathbb{S}^1$ and $T_{(q,\alpha)}(\overline{M}\times \mathbb{S}^1)\equiv T_q \overline{M}\oplus T_\alpha \mathbb{S}^1$, $q\in \overline{M}$, $\alpha\in \mathbb{S}^1$, for each tangent vector $v\in T_pS$ (or normal vector $w\in T_pS^\perp$) there is a canonical tangent vector $\hat{v}=(v,0)\in T_{(p,\alpha)}(S\times \mathbb{S}^ 1)$ (or $\hat{w}=(w,0)\in T_{(p,\alpha)}(S\times \mathbb{S}^1)^ \perp$). 
Moreover, it is clear that if  $S$ is minimal in $\overline{M}$, then $S\times \mathbb{S}^1$ is minimal in $\overline{M}\times \mathbb{S}^1$.

Finally, note that a similar procedure can be made is the submanifold is a geodesic.

} 
\end{rema}
We are now in position to state the first characterization result,

\begin{teor} \label{T1}
In a monotone Riemannian manifold $(I\times F,\beta \, dt^2 + g_{_t})$, every compact minimal submanifold
must be contained in a level hypersurface $t= const$.
Moreover, in the case of codimension higher than one, $S$ is a minimal submanifold of $(F,g_{_{t_0}})$, for some
$t_{_0}\in I$.
\end{teor}

\noindent\emph{Proof.}
Let $(S,g)$ be an $n(\geq 3)$-dimensional minimal submanifold under the assumptions. Endow $S$
with the conformal metric $\widetilde{g}= \beta^{(n-2)/2} g$. The function $\tau$ on $S$ satisfies equation
(\ref{laplacianoconforme}). We will see that this function is $\widetilde{g}$-superharmonic (if the Riemannian manifold is non-shrinking) or $\widetilde{g}$-subharmonic
(if it is assumed the non-expanding hypothesis). For this purpose, note that it is enough to see that $\sum_i \partial_t g_{_{ii}} - (\partial_t g_{_t})(N_i^F, N_i^F)$ is 
non-negative (resp. non-positive) if $I\times F$ is non-shrinking (resp. non-expanding) throughout $\partial_t$. In order
to prove that, consider on the Riemannian manifold the tensor field $\xi$ defined at each point by $\xi(u,v)=(\partial_t g_{_t})(d\pi_{_F}(u),d\pi_{_F}(v))$, for all tangent
vectors $u,v$.  Hence, for $p\in S$, 
$ \sum_i \partial_t g_{_{ii}} - (\partial_t g_{_t})(N_i^F, N_i^F) = \mathrm{tr}(\xi|_{T_p S})$. Now, it is clear that, from the non-shrinking assumption, we have that this term is non-negative (non-positive with the other kind of hypothesis).

Now, suppose $S$ has dimension at most $2$. Then, with the same procedure stated in Remark \ref{dimensionmas}, we realize that 
the conformal change can be applied. We conclude that $\tau$ must be constant.

Finally, assume $S$ is contained in a hypersurface $t=t_{_0}$, $t_{_0}\in I$. Taking into account equation (\ref{laplacianoconforme}) and
the tensor $\xi$, we arrive to
$$
\left( \partial_t g_{_t}\right)_{{t_{_0}}} (X,X)=0 \, ,
$$ for any $X \in \mathfrak{X} (S)$ (observe here that $\mathfrak{X} (S) \subset \mathfrak{X}(F)$). The 
Koszul formula can be called to obtain that $S$ is also a minimal submanifold of $(F,g_{_{t_0}})$. \hfill{$\Box$}

\begin{rema} \label{nota1}
{\rm
 Previous result is a complete classification of compact minimal hypersurfaces. Moreover, in that case, the hypersurface
must be totally geodesic. It follows as an application of the following general formula
$$
\overline{g}\left( \overline{\nabla}_X \frac{\partial_t}{\sqrt{\beta}}, Y  \right) = \frac{1}{2\sqrt{\beta}} (\partial_t g_{_t})(X,Y) \, ,
$$ for all $X,Y$ tangent vectors to $F$. 
Clearly, in order to exist such a hypersurface, $F$ must be compact.}
\end{rema}

\begin{rema}\label{noexis}
{\rm The proof of previous result allows us to
have some criteria to decide if a hypersurface $t=t_{_1}$ can contain a mininal submanifold. In fact, a necessary condition for the hypersurface
$t=t_{_1}$ to contain a minimal submanifold of dimension $n$ is that the tensor field
$(\partial_t g_{_t})|_{_{t_0}}$, at any point of the hypersurface $t=t_{_1}$, is degenerate with dimension of its radical at least $n$. The following
consequence also points in this direction.
 }
\end{rema}

\

Let us assume that $(M,g)$ admits a global descomposition as monotone throughtout two vector fields. This means that
$(M,g)$ is isometric to $(I_1\times F,dt_1^2+g_{t_1})$ and $(I_2\times F,dt_2^2+g_{t_2})$, $I_i\subseteq \mathbb{R}$, for $i=1,2$. If 
$g(\partial_{t_1},\partial_{t_2})^2\neq g(\partial_{t_1},\partial_{t_1}) \, g(\partial_{t_2},\partial_{t_2})$,
then $(M,g)$ is monotone throughout 2 non-collinear vector fields. It is trivial to extend to an arbitrary number of non-collinear vector fields.

\begin{coro}
Let $(M^n,\overline{g})$ be a Riemannian manifold which is monotone throughout $q$ ($\leq n$) non-collinear vector fields. Suppose that
for any point $p\in M$, the tangent vectors at $p$ generate a $q$-dimensional subspace of $T_p M$. There exists no
compact minimal submanifold of dimension at most $n-p+1$.
\end{coro}

\vspace{2mm}

On one hand, recalling Remark \ref{noexis}, if that condition does not hold for any hypersurface $t=t_{_0}$, then
Theorem \ref{T1} reads as a non-existence result. The strictly monotones Riemannian manifolds have $(\partial_t g_{_t})$
definite positive or negative. In particular, no hypersurface $t=const.$ can contain a minimal submanifold.

\begin{coro}\label{strictmonotone}
In a strictly monotone Riemannian manifold $(I\times F,\overline{g}=\beta dt^2+g_{_t})$ exist no compact minimal submanifolds.
\end{coro}

On the other hand, we deepen the last conclusion of Theorem \ref{T1}.  Consider a $2$-parametric Riemannian metrics on a 
manifold, $(F,g_{_{s,t}})$, $s\in J$, $t\in I$, where $I$ and $J$ are open intervals of the real line (perhaps the whole real line).
Then, given two function $\beta,\gamma \in C^\infty(I\times J \times F)$, we can build the Riemannian manifold
$(I\times J \times F,\beta dt^2 +\gamma ds^2 +g_{_{s,t}})$. Assuming this Riemannian manifold (that we may agree to write as $(I\times (J\times F),\beta dt^2+g_{_s,t})$) 
is under the assumptions of Theorem \ref{T1}, we arrive to compact minimal submanifolds (higher codimension than one) to 
be contained in an hypersurface $t=const$. We find that $S$ 
is minimal in the Riemannian manifold $(J\times F,\gamma|_{t_0} ds^2+g_{_{t=t_0,s}})$. But again
Theorem \ref{T1} can be used to state that it must be contained in a submanifold
$s=s_0$, for certain $s_0 \in J$. Observe that this process can be iterated indefinitely.

We define the following class
of Riemannian manifolds. Consider $m$ intervals of the real line $I_{_i}$, $i=1,\ldots, m$, with
a coordinate atlas $t_{_i} \in I_{_i}$. Consider a Riemannian manifold whose Riemannian metric depends on $m$ parameters, $(F,g_{_{t_{_1},\ldots, t_{_m}}})$, where 
$t_{_i}\in I_{_i}$ for $i=1,\ldots m$. Take also a ordered set of $m$ functions $\beta_{_i}\in C^\infty (I_{_1}\times \ldots \times I_{_m}\times F)$.
With these ingredients we can build the Riemannian manifold 
$(I_{_1}\times \ldots \times I_{_m}\times F,\beta_{_1} dt_{_1}^2+\ldots +\beta_{_m} dt_{_m}^2+g_{_{t_{_1},\ldots,t_{_m}}})$. For this class of Riemannian
manifolds, we have

\begin{coro} \label{moreover}
Let $(I_{_1}\times \ldots \times I_{_m} \times F, \overline{g}=\beta_{_1} dt_{_1}^2+ \beta_{_m} dt_{_m}^2 +g_{_{t_1,\ldots,t_m}})$ be a Riemannian manifold.
Assume all $\mathcal{L}_{\partial_{t_i}} \overline{g}$, $i=1,\ldots, m$, are together definite non-negative or non-positive tensor fields.

Let $B\subset \left\{1,2,\ldots, n\right\}$ be, with $m$ elements. Then, the only compact minimal submanifolds of codimension $m$ must be contained in a submanifold $t_{_i}=const.:i \in B$.

Moreover, if all $\mathcal{L}_{\partial_{t_i}} \overline{g}$, $i=1,\ldots, m$ are definite positive or negative $2$-covariant tensor fields, then it does not exist
compact minimal submanifolds.
\end{coro}   

The previous corollary can be specialized to the case where $F$ is an interval of the real line.
\begin{coro}
Let $(M^n,\overline{g})$ be a Riemannian manifold such that it is isometric to $(\Pi_{i=1}^n I_i, \overline{g}=\sum_{i=1}^n f_i \, dx_i^2)$, where $I_i\subseteq \mathbb{R}$ and $f_i\in C^\infty(\Pi_{i=1}^n I_i)$. Assume that for any $i,j\in 1, \ldots, n$, $\partial_{x_i} f_j$ is a monotonic function (resp. strictly monotonic
function). 

Let $B\subset \left\{1,\ldots, n\right\}$ be, with  $m$ element, then the compact minimal submanifolds of codimension $m$ are of the form  $\left\{x_l=const.: l\in B
\right\}$  such that
$\partial_{x_l}f_i|_{B}=0$ for any $i$ (resp. Then there exists no compact minimal submanifold).
\end{coro}

\vspace{2mm}

Now, we focus on Riemannian manifolds of constant sectional curvature.
Recall that any simply-connected Riemannian manifold of constant sectional curvature is, removing some points, isometric
to: $\mathbb{R}^n-p=((0,\infty)\times \mathbb{S}^n,dr^2+ r^2 \, g_{_{\mathbb{S}^{n-1}}})$, $\mathbb{H}^n(-k)-q=((0,\infty)\times \mathbb{S}^n,dr^2+\sqrt{k}^{-1} \cosh^2 (\sqrt{k} \, r) \, g_{_{\mathbb{S}^{n-1}}})$
or $\mathbb{S}^n(k)-\left\{n,s\right\}=((0,\frac{1}{\sqrt{k}})\times \mathbb{S}^n,dr^2+ \sqrt{k}^{-1} \sin^2 (\sqrt{k} r) \, g_{_{\mathbb{S}^{n-1}}})$ (here $ g_{_{\mathbb{S}^{n-1}}}$ denotes the 
canonical metric of the round sphere of radius $1$). In the non-compact case, assume there exists a compact minimal submanifold. Then, taking
$p$ not belonging to the submanifold, we may apply our results to find a contradiction. In the compact case, $\mathbb{S}^n(k)-n$ is expanding
thoughout $\partial_r$ until some value of $r$. However, we can assume that we restrict on a part of the total Riemannian manifold. Denote by
$B_p(r)$ the open geodesic ball centered at $p$ with radius $r$. It is easy to see that $B_p(\mathrm{diam}(\mathbb{S}^n(k), g_{\mathbb{S}^n(k)})/2)-p$,
endowed with the restricted metric is expanding. Assume that 
there exists a compact minimal submanifold $S$ in $B_p(\mathrm{diam}(\mathbb{S}^n(k), g_{\mathbb{S}^n(k)})/2)$. If $p\notin S$, we may apply our previous results to obtain a contradiction.
If $p\in S$, we may find points $q_i$ and numbers $s_i$ ($\leq \mathrm{diam}(\mathrm{S^n(k)}, g_{\mathbb{S}^n(k)})/2$)  such that: \emph{i)} $q_i\notin S$, \emph{ii)} $B_q(s_i) \subset B_p(\mathrm{diam}(\mathbb{S}^n(k), g_{\mathbb{S}^n(k)})/2)$, and \emph{iii)} $S \subset B_q(s_i)$. Applying the results to this set of points $q_i$,
we find again a contradiction with the existence of such a compact minimal submanifold. We have proved,

\begin{coro} \label{curvaturaconstante}
No simply-connected Riemannian manifold with non-positive constant sectional curvature admits a compact minimal submanifold.
In a round sphere, there exists no compact minimal submanifold contained in an open ball of radius the half of its diameter.
\end{coro}

\begin{rema} \label{greatcircles} {\rm
\textbf{a)} For a round sphere, observe that
any geodesic sphere of radius the half of the diameter of the Riemannian manifold is a minimal hypersurface. These minimal submanifolds are nice counterexamples to see
that our kind of assumptions are needed.

\textbf{b)} On the other hand, some topological assumption is necessary, as the simply-connectedness. Consider the torus
$T^3$. It is clear that there exist compact minimal submanifolds.     
} 
\end{rema}

\subsection{Change in the monotonic behaviour}

In this subsection we are interested in the case in which the monotonicity of the expanding
behavior of a Riemannian manifold changes. We will require the existence of 
a $t_{_0}\in I$ which divides the manifold into two parts which has different behavior.

\begin{teor} \label{T2}
 Let $(I\times F,\overline{g}=\beta dt^2+g_{_t})$ be a Riemannian manifold. Assume there exists $t_0\in I$ such that the manifold is non-expanding in the region
$t\leq t_0$ and non-shrinking in the region $t \geq t_0$.

 The only compact minimal submanifolds must be contained in a level
 hypersurface $t=const$.
\end{teor}

\noindent\emph{Proof.}
First, observe that on a compact Riemannian manifold $(M,g)$, the only functions such that $f\Delta f \geq 0$ are the constant
functions. In fact, from $\Delta f^2 = 2 |\nabla f|^2 + 2f\Delta f $ we get that $f^2$ is superharmonic, and
therefore constant from the compactness of $(M,g)$.

Assume dimension of the minimal submanifold at least $3$. From equation (\ref{laplacianoconforme}), the function $\tau$ satisfies
$$
(\tau-t_0) \widetilde{\Delta} (\tau -t_0) \geq 0 \, .
$$ Then, $\tau$ must be constant. The $2$-dimensional case follows analogously using an extension argument as used in the proof of Theorem \ref{T1}.
\hfill{$\Box$}


\

Following Remark \ref{greatcircles}, great circles of round spheres $\mathbb{S}^n(k)$ are counterexamples when the behavior is not as stated. 
In fact, writing $\mathbb{S}^n(k)n-\left\{n,s\right\}$ as above, then it is non-shrinking in the region $r\leq \frac{1}{\sqrt{k}}$ and non-expanding in the
region $r\geq \frac{1}{\sqrt{k}}$.

\vspace{1mm}

To end this section, we provide an application to warped product Riemannian manifolds.
Consider an interval of the real line $(I,dt^2)$, a Riemannian manifold $(F,g_{_F})$, and a function $f$ on $I$. The warped product Riemannian manifold
is the product manifold $I \times F$ endowed with metric $dt^2+f(t)^2 \, g_{_F}$. Following \cite{ON}, we denote this manifold as $I\times_f F$. 

\begin{coro}
Let $I\times_f F$ be a warped product Riemannian manifold. Assume that $f(t)$ has not a local maximum value. The only compact minimal
submanifolds must be contained in a level hypersurface $t=t_{_0}$ such that $f'(t_{_0})=0$. Moreover (when the codimension is greater than one), they must be
minimal submanifolds of $F$.
\end{coro} 

\noindent\emph{Proof.}
If $f$ is monotone, it may be used Theorem \ref{T1}. Otherwise, Theorem \ref{T2} can be called.
\hfill{$\Box$}

\vspace{2mm}

Observe that if $f$ has not critical points, such that minimal submanifolds cannot exist. Previous result can be combined with
Corollary \ref{curvaturaconstante} to obtain,
\begin{coro}
Let $I\times_f F$ be a warped product Riemannian manifold such that $(F,g_{_F})$ is simply-connected, it has non-positive constant sectional 
curvature and $f$ does not attain a maximum value. There exists no compact minimal submanifold of codimension bigger than $1$.
\end{coro}

To close this section, we can relax the hypothesis making use of a future result, Corollary \ref{corCH},

\begin{coro}
Let $I\times_f F$ be a warped product Riemannian manifold such that $(F,g_{_F})$ is a complete Cartan-Hadamard manifold and $f$ does not attain a maximum value. There exists no simply-connected compact minimal submanifold of codimension bigger than $1$.
\end{coro}

\section{Controlling a volume function}

In this section, we focus on different assumptions on the Riemannian manifold $(I\times F, dt^2+g_{_t})$. The main difference between previous sections is that
there it was required a common global behavior of $\mathcal{L}_{\partial_t} \overline{g}$, while here we only require that this tensor field
is semi-definite at each point. Hence, we may consider here Riemannian manifolds for which previous results cannot apply.

We come back to equation (\ref{divt}), and follow the notation from there. 
In the analysis of  $\overline{\mathrm{div}}(\partial_t)$ we may have approached in a different way. Consider 
a coordinate system $(\mathcal{V},(t,x_1,\ldots , x_m))$ of $\overline{M}$ and take the
canonical Riemannian volume element $\Omega \in \Lambda^{n+1}(\overline{M})$
$$
\Omega =  \sqrt{\mathrm{det}(g_{_t}(\partial_{x_i},\partial_{x_j}))} \, dt \wedge dx_1 \wedge \ldots \wedge dx_m \, . 
$$ We can write
$$
\overline{\mathrm{div}} \left( \partial_t  \right) = \frac{1}{2} \partial_t \log (\mathrm{det}(\partial_{x_i},\partial_{x_j}))  \, .
$$ Hence, the function $\eta := \partial_t \log (\mathrm{det}(\partial_{x_i},\partial_{x_j}))$ is globally defined and independent of the choice of coordinates.

\ 

On a minimal submanifold $S$ in $(I\times F, dt^2+g_{_t})$, consider the vector field $Y:= \eta \, \partial_t^\top$. From (\ref{laplaciano}), it obeys

\begin{equation} \label{divY}
\mathrm{div} \left(  Y \right) = \partial_t^\top \eta + \eta \, \Delta \tau \, .
\end{equation} The acute angle function can help us to write
$$
\partial_t^\top = \sin^2 \theta \, \partial_t + \sin \theta \cos \theta\,  u \, ,
$$ where the vector field $u$ is unitary and satisfies $g(u,\partial_t)=0$. Now, equation (\ref{divY}) leads to
\begin{equation} \label{divYY}
\mathrm{div} \left(  Y \right) = \sin^2 \theta \left( \partial_t \eta + \cot \theta \, u\left(\eta \right)  \right) + \eta \Delta \tau \, .
\end{equation}

\begin{teor}
Let $(I\times F, dt^2+g_{_t})$ be a Riemannian manifold such that, at each point, $\partial_t g_{_t}$ is semi-definite. Assume $\eta$
satisfies $\partial_t \eta \geq \sigma |\nabla^F \eta|$, for some $\sigma \in \mathbb{R}^+$. 
There exists no compact minimal submanifold with
acute angle function satisfying $\tan \theta \geq \sigma^{-1}$.
\end{teor}

\noindent\emph{Proof.}
First, we will see that $\mathrm{div} \left(  Y \right) \geq 0$. From (\ref{divYY}), it is enough to show that 
$\eta \Delta \tau \geq 0$. Taking into account that $\partial_t g_t$ semi-definite, it is easy to obtain the assertion.

Hence, the Divergence Theorem leads to $\eta =0$, since $\theta \geq \epsilon > 0$, for some positive constant $\epsilon$. Then, from
(\ref{laplaciano}), it is found that $\Delta \tau=0$, so $\tau$ must be constant. Contradiction.
\hfill{$\Box$}

\vspace{2mm}

Previous result may be geometrically interpreted as an impossibility to build minimal submanifolds. More precisely, it cannot be exhibited
any minimal submanifold contained in an hypersurface whose acute angle satisfies $\tan \theta \geq \sigma^{-1}$. Observe that the lower
estimation depends only on the geometry of the ambient Riemannian manifold.

\section{Uniqueness results in Riemannian manifolds with symmetries}

Let $(M, \overline{g})$ be an $(n+1)$-dimensional Riemannian manifold which possesses a Killing vector field $K$. 
If that vector field fulfill some assumptions, then we get a topological and geometrical description of Riemannian manifold.
This fact can be consulted in \cite[Proposition 1]{RomeroRubio}. Since the length of the proof, we reproduce the arguments here.

The Frobenius theorem asserts that
the orthogonal distribution of $K$ is integrable if and only if $K$ is irrotational. Locally, if $\Sigma$ is an open set of an integral leaf of $K^\bot$, $P$,
then $M$ is locally isometric to the product of $(\Sigma,g_{_\Sigma})$ with $(I,dt^2)$ endowed with metric $h^2 dt^2+g_{_\Sigma}$, where 
$h\in C^\infty (\Sigma)$. It is not difficult to see that the vector field $\frac{1}{h^2} K$ is locally a gradient vector field. Furthermore,
assuming that $M$ is simply-connected, then $\frac{1}{h^2} K$ is globally a gradient, $\mathrm{grad} \, l = \frac{1}{h^2}K$, for certain
$l\in C^\infty (M)$. Observe that
the metrically equivalent $1$-form $w$ associated to the vector field $\frac{1}{h^2}K$ is exact, $dl=w$. 

Let's denote by $\phi (t,p)$ the global flow of $K$. Then $\frac{d}{dt}l(\phi(t,p))=1$. Thus, the integral curves of $K$ cross each leaf
of $K^\bot$ only one time. We have that the map
$$
\varphi: P\times \mathbb{R} \rightarrow M, \quad \varphi(p,t)=\phi(t,p) 
$$ is an isometry. We have arrived to \cite[Proposition 1]{RomeroRubio}
\begin{prop}
Let $M$ be a Riemannian manifold which admits an irrational nowhere zero Killing vector field $K$. If $M$ is simply-connected and $K$ is complete,
then $M$ is globally isometric to a warped product $P\times_h \mathbb{R}$, where $P$ is a leaf of the foliation $K^\bot$ and $h=|K|$.
\end{prop} Recalling Theorem \ref{T1}, previous proposition leads to the following result, which generalizes \cite[Theorem 3 and 4]{RomeroRubio}),

\begin{teor} \label{TSimetrias}
Let $M$ be a complete simply-connected Riemannian manifold which admits an irrotational nowhere-zero Killing vector field $K$. 

Every minimal compact submanifold must be contained in a leaf of the foliation $K^\bot$.
\end{teor} 

Consider $x:S\rightarrow M$ be an immersion of $S$ in $(M,g)$. If $S$ is simply-connected and compact, then, in the universal
Riemannian covering of $(M,g)$, $(\widetilde{M},\widetilde{g})$ (take $\pi:\widetilde{M}\rightarrow M$ a covering map), we 
have a unique immersion $\widetilde{x}:S\rightarrow \widetilde{M}$ such that $\widetilde{x}\circ \pi = x$. Note that $x:S\rightarrow M$
is minimal if and only if $\widetilde{x}:S \rightarrow \widetilde{M}$ is so. Hence, the simply-connected assumption can be assumed on the compact minimal submanifold.

\begin{coro}
Let $M$ be a complete Riemannian manifold which admits an irrotational nowhere-zero Killing vector field $K$. 

Every simply-connected minimal compact submanifold must be contained in a leaf of the foliation $K^\bot$.

\end{coro}

To end this section, we desire remark that other kind of splitting theorems (see, for instance, \cite{Ponge}) can be equally 
combined with our results in order to obtain other uniqueness results.

\section{Applications to Geometric Analysis}

In this section, we study the case in which the minimal submanifold is a graph on $F$. Several considerations lead us to
state some uniqueness results for certain families of PDEs.
First, in order to be used latter, we present a technical lemma,

\begin{lema}\label{ultimolema}
Let $ (I\times F,\beta dt^2+g_{_t})$ be a Riemannian manifold and $S$ a compact hypersurface. For each function 
$h \in C^\infty (S)$, there exists a function in the ambient manifold, $\alpha$,
such that
\begin{equation} \label{nablaalfa}
 \nabla \alpha |_{S}^{\bot} = h N \, .
\end{equation}
\end{lema}

\noindent\emph{Proof.}
For each $p\in S$, let $\gamma_p (s)$, $s\in J$, the unique geodesic which satisfies $\gamma_p (0) =p$ and $\gamma'_p (s)=N(p)$. 
Consider the tubular neighbourhood of $S$, $U= \left\{ \gamma_p (s) : s \in J, p \in S \right\}$. The flow
associated to the geodesics $\gamma_p (s)$ on $J\times S$ is given by $\phi (s,p) = \gamma_p (s)$, where $\phi$ is
bijective. In $U$ we define the function
$\alpha (\phi^{-1}(s,p))$, $(s,p)\in J\times S$, by $\alpha (\phi^{-1}(s,p))=f(p)$, that is, $\alpha$ is constant along the geodesics $\gamma_p (s)$
and on $S$ it coincides with $f$. The normal gradient satisfies (\ref{nablaalfa}), since $\overline{g}(N,\nabla \alpha)=0$.

Now, let $\xi$ be a function on $I\times F$ such that $0\leq \phi(p)\leq 1$, for all $p\in I\times F$, and which satisfies (see Corollary in Section 1.11 of \cite{W}),

\textit{i)} $\xi(p)=1$ if $p\in \{\gamma_t(p): t\in J',\, p\in S\}$, being $J'\subset J$ an closed interval with $0\in J'$.

\textit{ii)} ${\rm supp}\, \xi\subset U$.

\vspace{2mm}
The function $\xi$ can be employed to extend $\alpha$ on all $I\times F$.
\hfill{$\Box$}

\vspace{2mm}

Consider an immersion $x:S^n \rightarrow (I\times F, \overline{g}=\beta dt^2+g_{_t})$. We can consider also the same immersion
when the ambient manifold is endowed with certain pointwise conformal metric, $\widetilde{x}:S\rightarrow (I\times F,\widetilde{g}=e^{2\alpha}\overline{g})$,
where $\alpha \in C^\infty (I\times F)$. The normal vector field of $S$ in $(I\times F,\overline{g})$, $N$, is related
with the same in $(I\times F,\widetilde{g})$, $\widetilde{N}$, by $\widetilde{N}=e^\alpha N$. Taking $\left\{E_i\right\}_{i=1}^n$ an orthonormal basis in $T_p S$,
$p\in S \subset (I\times F, \overline{g})$, then $\left\{e^\alpha E_i\right\}=\left\{\widetilde{E}_i\right\}$, $i=1,\ldots, n$, is an orthonormal basis in $T_p S$,
$p \in S\subset (I\times F,\widetilde{g})$. Denoting by $H$ and $\widetilde{H}$ the mean curvature function of $S$ in $(I\times F,\overline{g})$ and 
in $(I\times F,\widetilde{g})$, respectively, it is found that
$$
n\widetilde{H}=\sum_{i=1}^{n} \widetilde{g}(\widetilde{\nabla}_{\widetilde{E}_i} \widetilde{N}, \widetilde{E}_i ) = \sum_{i=1}^{n} e^{-\alpha} \overline{g}(\widetilde{\nabla}_{E_i} N,E_i) \, ,
$$ where $\widetilde{\nabla}$ is the Levi-Civita connection of $\widetilde{g}$. If $\overline{\nabla}$ denotes the Levi-Civita connection of $\overline{g}$,
from previous equation it follows
$$
e^{\alpha} \widetilde{H} = H + \overline{g}(\overline{\nabla} \alpha, N) \, .
$$ 

\begin{rema} \label{multipleexistencia}
{\rm
From Lemma \ref{ultimolema}, we are able to build compact minimal hypersurfaces. In fact, suppose given a compact minimal hypersurface: $x:S\rightarrow M$
in a Riemannian manifold $(M,\overline{g})$. Consider the mean curvature function on $S$. Then, there exists $\alpha \in C^\infty(M)$ such that $x:S\rightarrow (M,e^{2\alpha} \, \overline{g})$ is a minimal hypersurface. 
}
\end{rema}

\ 

We can apply this conformal change when $\alpha$ does not depend on the $t$-coordinate. 
Note that, in this case, the conformal change does not affect the expanding or contractive behaviour of the ambient
Riemannian manifold. 
Any function 
$u\in C^\infty (F)$ defines a  graph $\Sigma_u$ on $F$ by $\Sigma_u = \left\{ (u(p),p)\in I\times F : p \in F \right\}$. 
Denote by $H(u)$ the mean curvature operator associated to $\Sigma_u$, and $N^F$ the projection onto $F$ of the normal vector field associated to its graph.

As an application of Theorem \ref{T1} and \ref{T2}, we give,

\begin{teor}
Let $(I\times F, \beta dt^2+g_{_t})$, $F$ compact, be a Riemannian manifold such that it is non-expanding or non-shrinking throguhtout $\partial_t$,
or there exists $t^* \in I$ such that the manifold is non-expanding in the region $t\leq t^*$ and non-shrinking in $t \geq t^*$.
Then, on $F$, the equation
\begin{equation} \label{operadorH}
H(u)=-g_{_F} (N_{_F},D\alpha ) \, ,
\end{equation} where $\alpha \in C^\infty (F)$ has no solutions unless $u$ is the constant function.
\end{teor}

\noindent\emph{Proof.}
Assume $u$ is a non-constant function obeying (\ref{operadorH}). Employ Lemma \ref{ultimolema} in order to extend $\alpha$
to a function in the ambient space. Then, $u$ defines a compact hypersurface in the manifold
$(I\times F, \overline{g}=\beta dt^2 +g_{_t})$. In the conformal manifold, $(I\times F, e^{2\alpha} \overline{g})$, $u$
is a compact minimal hypersurface. The proof ends noting that we are now in position to use Theorem \ref{T1} or \ref{T2}.
\hfill{$\Box$}

\vspace{2mm}

Note that the associated Bernstein type problems appear when $\alpha = 0$ in previous result.

To put some concrete PDE for which previous theorem apply, we will compute explicitely the expression of $H(u)$ in some well-known relevant cases. Before that,
we need to develop a general equation.
Let $(M,\overline{g}=\beta dt^2+ g_{_t})$ be a Riemannian
manifold. Each function $u \in C^\infty (F)$ defines a graph $\Sigma_{_u}=\left\{ (u(p),p)\in I\times F\right\}$. For any function $f\in C^\infty (I\times F)$,
let's denote by $\nabla^F f$ the gradient on $F$, that is, $\nabla^F f = \overline{\nabla} f+\overline{g}(f,\frac{1}{\sqrt{\beta}}\frac{\partial_t}{\sqrt{\beta}})$.
Then, the normal vector field of a graph $\Sigma_{_u}$ on $S_{_0}$ is given by
$$
N= \frac{1}{\sqrt{\frac{1}{\beta}+|\nabla^F u|^2}} \left\{ \nabla^F u +\frac{1}{\beta} \partial_t   \right\} \, .
$$ The mean curvature of $\Sigma_{_u}$ is not difficult to compute,
\begin{eqnarray*}
nH &=& \overline{\mathrm{div}} \, N = \mathrm{div}_F (N^F) +\overline{g}\left( \overline{\nabla}_{\frac{\partial_t}{\sqrt{\beta}}} N^F, \frac{\partial_t}{\sqrt{\beta}}  \right) + \overline{\mathrm{div}} \left( \overline{g}(N,\frac{\partial_t}{\sqrt{\beta}}) \frac{\partial_t}{\sqrt{\beta}}   \right) \\
&=& \mathrm{div}_F \left( \frac{\nabla^F u}{\sqrt{\frac{1}{\beta} +|\nabla^F u|^2}}  \right) + 
\overline{g}\left( \frac{\nabla^F u}{\sqrt{\frac{1}{\beta} +|\nabla^F u|^2}},\frac{1}{2} \nabla^F \log \beta  \right)\\
&& + \overline{g}(N,\frac{\partial_t}{\sqrt{\beta}}) \frac{1}{\sqrt{\beta}} \partial_t \log \mathrm{vol}_{slice} \, .
\end{eqnarray*} Previous expression allows us to determine the operator $H(u)$ when  the explicit form of the metric of the ambient space
is known. 

\begin{ejem} \label{warpedproducts}
{\rm (see, for instance, \cite{RomeroRubio}). Let $I\times_f F^n$ be a warped product, where $f\in C^\infty (I)$. The minimal
hypersurface equation on $F$ is
$$
\mathrm{div} \left(  \frac{Du}{f(u) \sqrt{f(u)^2+|Du|^2}}   \right) = \frac{f'(u)}{\sqrt{f(u)^2+|Du|^2}}  \, \left\{ n- \frac{|Du|^2}{f(u)^2} \right\} \, .
$$
}
\end{ejem}



\begin{ejem} \label{ej21}
{\rm
Let $(I \times F, h^2 dt^2+g_{_F})$, where $h\in C^\infty(I\times F)$ is a positive function. Then, 
the minimal hypersurface equation on $F$ is given by
$$
\mathrm{div} \left( \frac{h \, D u}{\sqrt{1 +h^2 \, |D u|^2}}  \right) = 
-\frac{1}{ \sqrt{1 +h^2 \, |D u|^2}} \, \overline{g}\left( D u, D h  \right) \, ,
$$ where $D$ is the Levi-Civita connection of $(F,g_{_F})$.
}
\end{ejem}

\begin{ejem}{\rm
Let $(I\times F_{_1}^{n_1} \times F_{_2}^{n_2},dt^2+f_{_1}^2 g_{F_1}+f_{_2}^2 g_{F_2})$. Following previous considerations,
$\nabla^F u = \sum_{i=1}^2 \frac{1}{f_{_i}^2} D^{F_i} u$, where $D^{F_i}$ is the Levi-Civita connection of $(F,g_{_{F_i}})$, 
for $i=1,2$. Then, we have that the minimal hypersurface equation, on $F_{_1}\times F_{_2}$, is
$$
\sum_{i,j=1}^2 \mathrm{div}_{F_j}   \left(  \phi \frac{1}{f_{_i}^2} D^{F_i} u   \right) = - \phi \left\{ n_1 (\log f_{_1})'(u) +n_2 (\log f_{_2})'(u) \right\} \, ,
$$ where $$
\phi^{-1}= \sqrt{1+f_{_1}^2 |D^{F_1}u|^2 + f_{_2}^2 |D^{F_2}u|^2} \, .
$$ The extension to a finite family of Riemannian manifolds follows easily.
}
\end{ejem} 

%

\section{Dirichlet problems}

Let us consider the problem of finding a (piece of) minimal hypersurface $\Sigma$ in $(I\times F, \beta \, dt^2+g_{_t}$ under the constrain
$\partial \Sigma \subset \left\{t= t_0\right\}$, $t_0 \in I$.

\begin{teor} \label{diri1}
Let $($$I\times F, \beta \, dt^2+g_{_t}$$)$ be a is non-shrinking $($resp. non-expanding$)$ Riemannian manifold. The only orientable compact minimal submanifold $\Sigma$ 
such that $\partial \Sigma \subset \left\{t_0 \right\}$, $t_0 \in I$ and $\tau \geq t_0$ $($resp. $\tau \leq t_0$ $)$, is 
a $($piece of$)$ $\left\{t=t_0\right\}$.
\end{teor}

\noindent\emph{Proof.}
Endow $\Sigma$ with the conformal metric $\widetilde{g}$ as in Section 3. Apply the dimensional extension as in Remark 1, if necessary. Assume
the ambient manifold is non-shrinking. The vector field $(\tau-t_{0})\widetilde{\nabla} \tau$ vanishes on $\partial \Sigma$. Using the Divergence Theorem,
$$
0 = \int_{\Sigma} \left( (\tau - t_0) \widetilde{\Delta} \tau + |\widetilde{\nabla} \tau|_{\widetilde{g}}^2 \right) \, \, d\widetilde{\Sigma} \, ,
$$ where $d\widetilde{\Sigma}$ denotes the area element of $(\Sigma,\widetilde{g})$. Since $(\tau - t_0) \widetilde{\Delta} \tau \geq 0$ from
hypothesis, we get that $\tau$ must be constant.
The non-increasing case follows analogously taking into account now the vector field $(t_0- \tau)\widetilde{\nabla} \tau$.
\hfill{$\Box$}

\ 

The previous theorem can be combined with several families of PDEs in order to produce uniqueness of Dirichlet problems. For instance,
from Example \ref{warpedproducts},

\begin{ejem} {\rm
Let $\Sigma$ be a compact domain of a Riemannian manifold $F$, with $\partial \Sigma \neq \emptyset$, and let $f:\mathbb{R} \rightarrow \mathbb{R}^+$
be a smooth non-decreasing (resp. non-increasing function). The only solution $u\in C^\infty (\Sigma)$ to

\begin{eqnarray*}
 \mathrm{div} \left(  \frac{Du}{f(u) \sqrt{f(u)^2+|Du|^2}}   \right) = \frac{f'(u)}{\sqrt{f(u)^2+|Du|^2}}  \, \left\{ n- \frac{|Du|^2}{f(u)^2} \right\}  \\
u \geq t_0  \quad ( \mathrm{resp.} \, \, u \leq t_0 \, ) \\
u(\partial \Sigma) = t_0 \, , \, 
\end{eqnarray*}

is $u=t_0$ if $f'(t_0)=0$. Otherwise, there is no solution.}
\end{ejem}

Using again the vector field $(\tau-t_{0})\widetilde{\nabla} \tau$ we can provide the following result,

\begin{teor}
Assume $($$I\times F, \beta \, dt^2+g_{_t}$$)$ is non-shrinking in $t\geq t_0$ and non-expanding in $t\leq t_0$. The only minimal hypersurface $\Sigma$ 
such that $\partial \Sigma \subset \left\{t_0 \right\}$, $t_0 \in I$ is 
a $($piece of$)$ $\left\{t=t_0\right\}$.
\end{teor}

We can particularize to the case in wich $\partial_t$ is a Killing vector field. We may use Example
\ref{ej21} to get the minimal hypersurface equation.
\begin{coro}
Let $\Sigma$ be a compact domain with boundary and $h\in C^\infty (\Sigma)$. The only solutions $u\in C^\infty(\Sigma)$ to
$$
\mathrm{div} \left( \frac{h \, D u}{\sqrt{1 +h^2 \, |D u|^2}}  \right) = 
-\frac{1}{ \sqrt{1 +h^2 \, |D u|^2}} \, \overline{g}\left( D u, D h  \right) \, ,
$$ with the Dirichlet boundary condition
$$
u = const. \quad \mathrm{on} \, \,  \partial \Sigma \, ,
$$ are the constant functions.

\end{coro}

\ 

Finally, we want to provide this corollary, which gives us geometrical information about the shape of minimal submanifolds,

\begin{coro} \label{nomaximos}
Assume $($ $I\times F, \beta \, dt^2+g_{_t}$ $)$ is a expanding (resp. contacting) Riemannian manifold throughout $\partial_t$. Let $S$
be a minimal hypersurface. Then, the function $\tau$ on $S$ does not attain a strict maximum value (resp. strict minimum value).
\end{coro}

\noindent\emph{Proof.}
Assume there exists $S$ not satisfying our conclusion, in the expanding case. Let us say that $\tau_0$ is a maximum value of $\tau$ at $p\in S$. We have that, for
certain $\delta>0$ small enough, there exists a simply-connected compact oriented subset $\Sigma$ of $S$ containing $p$ and whose boundary lies in ${\tau_0-\delta}$. Inside this subset,
$\tau\geq t_0-\delta$. Apply Theorem \ref{diri1} in order to get a contradiction. Similar arguments can be applied to prove the contracting case.
\hfill{$\Box$}

\section{Applications to Riemannian Geometry}

Let $(M^n,g)$ be a Riemannian manifold and consider $p\in M$. Take a normal neighbourhood
$\overline{U}$, with $p\in \overline{U}$, and normal coordinates $(x_1,\ldots , x_n)$ in $\overline{U}$.
Let $B(p,\epsilon)$ be a geodesic ball centered at $p$ with radius $\epsilon$ and such that $B(p,\epsilon)\subseteq \overline{U}$.
Denote by $\eta =\sqrt{g_{ij}(x)} \, dx_1 \wedge \ldots \wedge dx_n$  the Riemannian canonical volume element on $\overline{U}$.
Then
$$
\lim_{\epsilon \rightarrow 0} \int_{\overline{B}(p, \epsilon)} \sqrt{g_{ij}(x)} \, dx_1 \wedge \ldots \wedge dx_n =0 \, ,
$$ since, by definition
$$
\int_{\overline{B}(p,\epsilon)} \sqrt{g_{ij}(x)} \, dx_1 \wedge \ldots \wedge dx_n = \int_{\overline{B}_{\mathbb{R}^n}(0,\epsilon)}\sqrt{g_{ij}(x)} \, dx_1 \wedge \ldots \wedge dx_n \, ,
$$ and $\sqrt{g_{ij}(x)}$ is bounded in each compact set. From this fact, we deduce that the volume of the geodesic balls goes to zero
when the radius does.

\ 

Take now the normal unitary radial vector field $N$, given by $N=\frac{d}{dt}\left( \exp_p (t \, u)   \right)$, for unitary vectors
$u\in T_p \mathbb{R}^n$. 

The $(n-1)$-form $i_N \, \eta$ is the canonical volume element in each geodesic sphere of $\overline{U}$, which we will denote by $E(p,\epsilon)$. Using
the Stokes Theorem,
$$
\int_{E(p,\epsilon)} i_N \, \eta = \int_{\overline{B}(p,\epsilon)} d(i_N \eta) \, .
$$ Taking into account that $d(i_N \, \eta)=h(q) \, \eta$, where $h(q)$ is a bounded function by compacity we can
assert that the $(n-1)$-volume of the geodesic sphere $E(p,\epsilon)$ tends to zero when $\epsilon\rightarrow 0$.

\ 

On the other hand, taking $\epsilon>0$ small enough, there exists a diffeomorphism given by
$$
\exp_p : B_{\mathbb{R}^n}(0,\epsilon) \rightarrow E(p,\epsilon) \, .
$$ We can consider the map
$$
f:(0,\infty)\times \mathbb{S}^{n-1} \rightarrow E(p,\epsilon)-\left\{ p\right\} 
$$ defined by
$$
f(r,v)=\exp_p (rv) \, .
$$ Taking into account the pull-back of the metric of $M$ (on $E(p,\epsilon)-\left\{ p\right\}$), and the Gauss Lemma, we obtain
$$
f^{*}(g)=dr^2+h_{(r,v)} \, .
$$ This produces a $1$-parametric families of Riemannian metrics on $\mathbb{S}^{n-1}$, $(\mathbb{S}^{n-1},h_r)$.

\ 

Now, we will see that, given a vector $u\in T_q(\mathbb{S}^{n-1})$, $h_r(u,u)$ grows with $r$ in a neighbourhood $(0,\epsilon)$ of $r$.
Given a such a vector $u$, under the isometry previously defined, we can consider the vector $df_{(r,q)}(0,u)$, which will be tangent to
a geodesic sphere of radius $r$.
Consider in $T_p M \equiv \mathbb{R}^n$, the radial unitary vector field $\overline{N}(x):=x/||x||$, defined up the origin point. For
each $x$ (fixed), take the $2$-plane $\Pi(R,u)= \mathrm{span} \left\{ R(x),u  \right\}$, where $R(x)=r \, x/||x||$.

\ 

Taking the image of this $2$-plane by the map $\exp_p$ (where it will be defined), and taking its intersection with $B(p,\epsilon)$, we obtain
a surface $D(p,\epsilon)$ embedded in $B(p,\epsilon)$. Moreover, its intersection of the image of the plane with $E(p,\epsilon)$ is a curve $C(p,\epsilon)$ in
the geodesic sphere, where at any point, the velocity of the curve is equal to $u$.

\ 

Taking into account all this considerations, we have that
$$
\lim_{\epsilon\rightarrow 0} \mathrm{length}\left( C(p,\epsilon)  \right) =0 \, .
$$

\ 

Finally, and considering $\epsilon$ small enough, it is easy to see that
$$
\lim_{r\rightarrow 0} h_r (u,u)=0 \, .
$$ We have proved,

\begin{prop}\label{local}
Let $(M,g)$ be a Riemannian manifold. Then, locally, it is expanding throughout certain vector field, in an open subset up a point. 
More precisely, for each point $p\in M$, there exists a $\delta_p\in \mathbb{R}^+$ such that $(B(p,\epsilon)-p,g)$, is expanding
throughout the radial polar geodesic vector field centered at $p$, for some $\epsilon_p>0$.
\end{prop}

\ 

For each point of a Riemannian manifold $(M,g)$, there exists an open geodesic ball of radius $\delta_p$ in which
it is expanding. For each point $p$, denote by $\overline{\delta_p}$ the supremum of such radius. Then, we have a function
on $M$, $h(p)=\overline{\delta_p}$ which it is continuous and positive in $M$. If $M$ is compact, it must have a minimum
$\delta_0>0$. We get,

\begin{teor}\label{encompactas}
For any compact Riemannian manifold $(M^n,g)$ there exists $\delta>0$ such that in any open geodesic ball of radius $\delta$
there exists no compact minimal submanifold.
\end{teor}

\noindent\emph{Proof.} Suppose there exists a point $p$ such that for any ball $B(p,\delta)$, $\delta>0$ there exists
a minimal compact submanifold $S_{\delta}$ contained in it. If we take $\delta=\delta_0/2$, Corollary \ref{strictmonotone} allows us
to state that $p$ belongs to $S_{\delta_0/2}$. Using the Sard theorem, we can choose a point $q$ in the ball $B(p,\delta_0/2)$ arbitrarily close
to $p$ and with $q\notin S_{\delta_0/2}$, in such a way that the ball $B(q,\delta_0)$ includes the ball $B(p,\delta_0/2)$. This leads to
a contradiction again with Corolary \ref{strictmonotone}.
\hfill{$\Box$}

\ 

\begin{coro}
For any $3$-dimensional compact Riemannian manifold there exists $\delta>0$ such that there exists no compact minimal surface in any open geodesic ball
of radius $\delta$. Neither exists a sequence of closed geodesics whose legths tend to zero.
\end{coro}

\begin{rema}{\rm
As we have said in the Introduction, given a Riemannian manifold $(M,g)$, we may define $\overline{\delta}(M,g)$
as the greatest value among those which satisfies Theorem \ref{encompactas} divided by the diameter of $(M,g)$. Then, Corollary \ref{curvaturaconstante} shows 
that $\overline{\delta}(\mathbb{S}^n(k),g_{\mathbb{S}^n(k)})=1/2$.  A question that arises naturally: is the round sphere the only compact Riemannian manifold $(M,g)$ such that $\overline{\delta}(M,g)=1/2$?
}
\end{rema}

\ 

Now, we focus on the case in which the ambient Riemannian manifold $(M,g)$ is a simply-connected Cartan-Hadamard manifold. In this case, it 
is well known that, for any $p\in M$ the map $\exp_p$ is a global diffeomorphism. Let $v,w\in T_p M$, $|v| \, |w|>0$ two tangent
vectors. Define $\Pi(u,v)=\left\{ \exp_p \left( u \right), u \in \mathrm{Span}\left\{ v,w \right\}, u\neq 0 \right\}$. Note that $\Pi(u,v)$
is a totally geodesic hypersurface, when endowed with the induced metric, $g_{\Pi(u,v)}$. Morever, it has non-positive Gauss curvature. 
Denote by $\partial_r$ the radial polar vector field of $(M,g)$ centered at $p$. Then,
$(M-p,g)$ is expanding throughout $\partial_r$ if and only if $(\Pi(u,v), g_{\Pi(u,v)})$ is expanding throughout $\partial_r|_{\Pi(u,v)}$ for
any $u,v$.  
By an abuse of notation, we may represent $(\Pi(u,v), g_{\Pi(u,v)})$ as $\left( (0,\infty)\times \mathbb{S}^1, ds^2+f(s,\theta)^2 d\theta^2  \right)$,
where $\partial_s = \partial_r|_{\Pi(u,v)}$. The Gauss curvature is computed to be
$$
K= -\frac{\partial^2_s f(s,\theta)}{f} \, .
$$ Then, $K\leq 0$ implies that $\partial^2_s f(s,\theta)\geq 0$. Since $\partial_s f(s,\theta)>0$ in the set $(0,\epsilon)\times \mathbb{S}^1$, we have that in 
$(0,\epsilon)\times \mathbb{S}^1$,  $\partial_s f(s,\theta)>0$.

 We have proved,

\begin{teor}\label{CH}
Let $(M,g)$ be a simply-connected complete Cartan-Hadamard Riemannian manifold. Then, $M-p$ is expanding throughout the radial geodesic vector field centered at $p$.
\end{teor}

Using now Corollary \ref{strictmonotone}, we conclude (compare with \cite[Corollary 2]{atsuji})

\begin{teor} \label{TeoCH}
Let $(M,g)$ be a simply-connected complete Cartan-Hadamard Riemannian manifold. It admits no compact minimal Riemannian manifold.
\end{teor}

At this point, see Remark \ref{greatcircles} for 

Assuming topological assumptions on the compact minimal submanifold, 
\begin{coro} \label{corCH}
No simply-connected compact manifold can be minimally immersed in a complete Cartan-Hadamard manifold.
\end{coro}

\

Finally, from our study we can get geometrical information about the shape of a minimal submanifold. To fix ideas, let 
us given an immersion of a compact submanifold $S$ in the Euclidean space $\mathbb{R}^n$.  From 
Remark \ref{multipleexistencia}, there
exists a conformal metric for which $S$ is minimal (via the same immersion). In particular, it can be exhibited Riemannian manifolds
which possess a minimal submanifold whose graph (seen as its immersion in 
$\mathbb{R}^n$) attains an extremum point. This fact does not occur if the Riemannian manifold is Cartan-Hadamard. We are able to 
talk about (local) extremum points of a graph in a complete Cartan-Hadamard manifold $(M,b)$ since the following fact. As 
detailed in Section 5, given an hypersurface $S$ in $(M,g)$, via their corresponding universal Riemannian covering maps, we get
that the universal Riemannian covering of $S$, $\tilde{S}$, is determined uniquely by an immersion in $\mathbb{R}^n$ (not necesarily endowed with its flat metric). 
Then, $S$ has a strict extremum point if (locally) the graph of $\tilde{S}$ has a strict extremum point.

It can be observed that if the ambient manifold is the Euclidean space, then the classical maximum principle may be used. However,
there exists some other manifolds for which the conclusion is not achieved.

We desire to close this paper with the following theorem,

\begin{teor} \label{noextremos}
Let $x:S \rightarrow M^n$ be a minimal submanifold in a complete Cartan-Hadamard manifold. The corresponding minimal immersion in the
universal Riemannian covering, $\tilde{x}:\tilde{S} \rightarrow \mathbb{R}^n$ cannot have a strict extremum point.
\end{teor}

\end{document}